\newcommand{\MailFile}[1]{}
\documentclass[a4paper,12pt]{article}
  \newcommand{\thedate}{October 12, 2019}

  \usepackage{email}
  \MailFile{coset_posets_arxiv_version.tex}
  \DoNotMail{ams*,*.fd,ulem.sty,enumitem.sty,babel.sty,*.ldf,verbatim.sty,inputenc.sty}

  \usepackage[utf8]{inputenc}
  \usepackage[DinA4]{buxlayout}
  \usepackage{amssymb,amsmath}
  \usepackage[american]{babel}
  \usepackage[GlobalNumbered]{buxmath}
  \usepackage{enumitem}
  \usepackage[normalem]{ulem}

  \InputIfFileExists{basic_notation.tex}{}{}
  \newvariable{\scrc}{\text{scrc}}
  \newvariable{\NumGen}{s}
  \newvariable{\GenProb}{P}
  \newvariable{\TheSet}{X}
  \newvariable{\TheIndexSet}{A}
  \newvariable{\TheIndex}{\alpha}
  \newvariable{\AltIndex}{\beta}
  \newvariable{\TheSubset}{U}
  \newvariable{\AltSubset}{V}
  \newvariable{\TheFam}{\mathcal{U}}
  \newvariable{\AltFam}{\mathcal{V}}
  \newvariable{\TheGrFam}{\mathcal{K}}
  \newvariable{\AltGrFam}{\mathcal{L}}
  \newvariable{\CosetFam}{\operatorname{Cosets}}
  \newvariable{\CosNerve}{\mathcal{NC}}
  \newvariable{\CosSimpl}{\mathcal{XC}}
  \newvariable{\CosPoset}{\mathcal{PC}}
  \newvariable{\FinIndexFam}{\mathcal{K}_{\text{fi}}}
  \newvariable{\HerFam}{\mathcal{H}}
  \newvariable{\TheNormalizer}{N}

  \newvariable{\TheSimplex}{\sigma}
  \newvariable{\Simpl}{\mathcal{X}}
  \newvariable{\Nerve}{\mathcal{N}}
  \newvariable{\Poset}{\mathcal{P}}
  \newvariable{\TheGroup}{G}
  \newvariable{\TheSubgroup}{H}
  \newvariable{\TheSubSubgroup}{K}
  \newvariable{\TheComplement}{M}
  \newvariable{\MinNormalSubGr}{M}
  \newvariable{\MaxNormalSubGr}{M}
  \newvariable{\FixNormalSubGr}{\widetilde{N}}
  \newvariable{\LargeNormalSubGr}{L}
  \newvariable{\PairFam}{\mathcal{L}}
  \newvariable{\TheProjection}{\pi}
  \newvariable{\TheNerveMap}{q}
  \newvariable{\Suspension}{\Sigma}
  \newvariable{\TheRest}{S}
  \newcommand{\subgroup}{\leq}
  \newcommand{\propersubgroup}{<}
  \newcommand{\normalsubgr}{\trianglelefteq}
  \newcommand{\propernormalsubgr}{\triangleleft}
  \newvariable{\TheGroupSet}{M}
  \newvariable{\TheGenSet}{X}
  \newvariable{\TheRelSet}{R}
  \newvariable{\TheGenerator}{a}
  \newvariable{\Sl}{\operatorname{SL}}
  \newvariable{\TheSetElement}{v}
  \newvariable{\TheGroupElement}{g}
  \newvariable{\TheLength}{k}
  \newvariable{\TheGraph}{\Gamma}
  \newcommand{\crossprod}{\times}
  \newvariable{\AltSetElement}{u}
  \global\let\Out\undefined
  \newvariable{\Out}{\operatorname{Out}}
  \global\let\Aut\udnefined
  \newvariable{\Aut}{\operatorname{Aut}}
  \newvariable{\FreeGroup}{F}
  \newvariable{\FinGroup}{Q}
  \newvariable{\TheDim}{d}
  \newvariable{\TheNumber}{m}
  \newvariable{\semidirectprod}{\rtimes}
  \newvariable{\PrimeNumbers}{\PPP}
  \newvariable{\IntNumbers}{\ZZZ}
  \newvariable{\RealNumbers}{\RRR}
  \newcommand{\freeprod}{*}
  \newvariable{\TheKernel}{K}
  \newvariable{\TheNormalSubGr}{N}
  \newvariable{\AltSubgroup}{M}
  \newvariable{\TheMaxSubgroup}{M}
  \newvariable{\Perm}{\operatorname{Perm}}
  \newvariable{\ThePrime}{p}
  \newvariable{\AltGroup}{\widetilde{G}}
  \newvariable{\TheProj}{\pi}
  \newvariable{\DiederGr}{D}
  \newcommand{\homotopic}{\simeq}
  \newcommand{\isomorphic}{\cong}
  \newcommand{\intersection}{\cap}
  \newvariable{\TheSimpleGr}{S}
  \newcommand{\subcomplex}{\leq}
  \newvariable{\TheFratSubGr}{\Phi}
  \newvariable{\NorFratSubGr}{\widetilde{\Phi}}
  \newvariable{\TheSubComplex}{Y}
  \newcommand{\pref}[1]{(\ref{#1})}
  \newvariable{\NumComplements}{c}
  
  \newvariable{\join}{*}
  \newvariable{\QuotGroup}{\overline{\TheGroup}}
  \newcommand{\isom}{\cong}

  \renewcommand{\notion}[1]{\emph{#1}}
  \newcommand{\bottomtext}[1]{%
    \begingroup
    \renewcommand{\thefootnote}{\relax}%
    \renewcommand{\footnotemark}{\relax}%
    \footnotetext{#1}
    \endgroup
  }

  \newcommand{\msc}[1]{\bottomtext{2010 Mathematics Subject Classification: #1.}}

\begin{document}
  \title{Coset Posets of Infinite Groups}
  \author{Kai-Uwe~Bux \and Cora~Welsch}
  \date{\thedate}
  \maketitle
  \begin{abstract}\noindent
    We consider the coset poset associated with the families of
    proper subgroups, proper subgroups of finite index, and
    proper normal subgroups of finite index. We investigate
    under which conditions those coset posets have contractible
    geometric realizations.
  \end{abstract}
  \msc{20E15 (20F32)}
  \bottomtext{Financial support by the {\small DFG} through
    the programs {\small SPP~2026} and {\small SFB~878} is gratefully
    acknowledged.}

  \noindent
  Let $\TheSet$ be a set and let
  \(
    \TheFam=\FamOf[\TheIndex\in\TheIndexSet]{\TheSubset[\TheIndex]}
  \)
  be a family of subsets. H.\,Abels and S.\,Holz associate three simplicial
  complexes to $\TheFam$, namely:
  \begin{itemize}\item[]
      \begin{itemize}
        \item[$\SimplOf{\TheFam}$] denotes the complex of finite subsets of
          $\TheSet$ contained within some $\TheSubset[\TheIndex]$.
        \item[$\NerveOf{\TheFam}$] is the nerve of the family $\TheFam$, i.e., a
          subset $\TheSimplex\subseteq\TheIndexSet$ is a simplex of
          $\NerveOf{\TheFam}$ provided that the intersection
          \(
          \TheSubset[\TheSimplex] := \Intersection[\TheIndex\in\TheSimplex]{
          \TheSubset[\TheIndex]
          }
          \)
          is nonempty.
        \item[$\PosetOf{\TheFam}$] is the geometric realization of the poset
          $\SetOf[{ \TheSubset[\TheIndex] }]{ \TheIndex\in\TheIndexSet }$
          where the order is given by inclusion.
          Thus, simplices of $\PosetOf{\TheFam}$ are finite $\subset$-chains in
          $\SetOf[{ \TheSubset[\TheIndex] }]{ \TheIndex\in\TheIndexSet }$.
      \end{itemize}
  \end{itemize}
  Abels--Holz show~\cite[Theorem~1.4]{AH}
  \begin{enumerate}
    \item[(a)]
      $\SimplOf{\TheFam}$ and $\NerveOf{\TheFam}$ are homotopy equivalent.
    \item[(b)]
      $\SimplOf{\TheFam}$, $\NerveOf{\TheFam}$, and $\PosetOf{\TheFam}$ are all
      homotopy equivalent provided that the family $\TheFam$ is closed with respect to
      non-empty intersection, i.e., if $\TheSubset,\AltSubset\in\SetOf[{
      \TheSubset[\TheIndex] }]{ \TheIndex\in\TheIndexSet }$ have non-empty
      intersection $\TheSubset\intersect\AltSubset\neq\emptyset$ then the intersection
      $\TheSubset\intersect\AltSubset$ also belongs to $\SetOf[{ \TheSubset[\TheIndex]
      }]{ \TheIndex\in\TheIndexSet }$.
  \end{enumerate}
  Note that $\SimplOf{\TheFam}$ does only depend on the set
  \(
    \SetOf[ {\TheSubset[\TheIndex]} ]{ \TheIndex\in\TheIndexSet }
  \).
  Hence the homotopy type of $\NerveOf{\TheFam}$ does not depend on the
  particular indexing of this collection of subsets. We shall therefore
  restrict ourselves to the case of self-indexing families where
  \(
    \TheIndex = \TheSubset[\TheIndex]
  \), i.e., $\TheIndexSet$ \emph{is} a set of subsets of $\TheSet$.

  Let $\TheGrFam$ be a set of subgroups of a fixed group $\TheGroup$ and
  let $\CosetFam[\TheGrFam]$ be the corresponding collection of cosets.
  We put:
  \[
    \CosSimplOf{\TheGrFam} := \SimplOf{ \CosetFam[\TheGrFam] }
    \qquad
    \CosNerveOf{\TheGrFam} := \NerveOf{ \CosetFam[\TheGrFam] }
    \qquad
    \CosPosetOf{\TheGrFam} := \PosetOf{ \CosetFam[\TheGrFam] }
  \]
  Specializing the above result of Abels--Holz to these families,
  one obtains the following consequence.
  \begin{fact}\label{ah}
    \begin{enumerate}
      \item[(a)]
        $\CosSimplOf{\TheGrFam}$ and $\CosNerveOf{\TheFam}$ are
        homotopy equivalent.
      \item[(b)]
        If $\TheGrFam$ is closed with respect to finite intersections,
        then $\CosSimplOf{\TheGrFam}$, $\CosNerveOf{\TheGrFam}$ and
        $\CosPosetOf{\TheGrFam}$ are homotopy equivalent.
    \end{enumerate}
  \end{fact}
  We shall be interested in the family $\HerFamOf{\TheGroup}$ of all proper
  subgroups in $\TheGroup$, the family $\HerFamOf[\text{fi}]{\TheGroup}$ of all
  its proper subgroups of finite index, and the family
  $\HerFamOf[\text{nor,fi}]{\TheGroup}$ of normal subgroups from
  $\HerFamOf[\text{fi}]{\TheGroup}$. All three families are closed with respect to
  finite intersections.
  We study the topology of the coset nerves
  $\CosNerveOf{\TheGroup} := \CosNerveOf{\HerFamOf{\TheGroup}}$ and
  $\CosNerveOf[\text{fi}]{\TheGroup} :=
  \CosNerveOf{\HerFamOf[\text{fi}]{\TheGroup}}$ and
  $\CosNerveOf[\text{nor,fi}]{\TheGroup} :=
  \CosNerveOf{\HerFamOf[\text{nor,fi}]{\TheGroup}}$.

  In particular, we
  shall investigate for which groups these nerves are contractible.
  For the families $\HerFamOf[\text{fi}]{\TheGroup}$ and
  $\HerFamOf[\text{nor,fi}]{\TheGroup}$ we have a satisfying answer to
  this problem.
  \begin{NewTh}<statement>{Theorem~\ref{equiv}}
    The coset nerve $\CosNerveOf[\text{fi}]{\TheGroup}$ is contractible
    if and only if the collection $\HerFamOf[\text{fi}]{\TheGroup}$ has
    infinitely many maximal elements.
  \end{NewTh}
  \begin{NewTh}<statement>{Theorem~\ref{normal-equiv}}
    The coset nerve $\CosNerveOf[\text{nor,fi}]{\TheGroup}$ is contractible
    if and only if the collection $\HerFamOf[\text{nor,fi}]{\TheGroup}$ has
    infinitely many maximal elements.
  \end{NewTh}
  For the coset nerve $\CosNerveOf{\TheGroup}$ associated to the family
  of all proper subgroups, we have only partial results, which we present
  in Section~\ref{all-subgroups}.

  It follows from Fact~\ref{ah} that
  the coset nerve $\CosNerveOf{\TheGroup}$ is homotopy equivalent
  to the \notion{coset poset}
  \(
    \CosPosetOf{\TheGroup} := \CosPosetOf{\HerFamOf{\TheGroup}}
  \)
  which is the set of all proper subgroups of $\TheGroup$ and their cosets,
  ordered by inclusion.

  The coset poset was introduced for finite groups by K.S.\,Brown in
  \cite{Brown}.  He considered the Euler characteristic of the coset poset since
  it is connected to the probabilistic zeta function.  The probabilistic zeta
  function is the reciprocal of the probability $\GenProbOf{\TheGroup,\NumGen}$
  that a randomly chosen ordered $\NumGen$-tuple from a finite group $\TheGroup$
  generates $\TheGroup$. Hall gave a formula for $\GenProbOf{\TheGroup,\NumGen}$
  as a finite Dirichlet series. In view of that formula, one can evaluate
  $\GenProbOf{\TheGroup,\NumGen}$ at an arbitrary complex number $\NumGen$.  Brown
  proved S.\,Bouc's observation that $\GenProbOf{\TheGroup,-1}$ can be interpreted
  as the negative reduced Euler characteristic of the coset poset
  $\CosPosetOf{\TheGroup}$.  This result motivated Brown to study the homotopy
  type of the coset poset, which (as Brown said) raised more questions than it
  answered.

  The question of simple connectivity was studied by D.A.\,Ramras in
  \cite{Ramras}.  Among other results, he proved that the coset poset
  $\CosPosetOf{\TheGroup}$ is contractible if $\TheGroup$ is an infinitely
  generated group.  The question of contractibility was answered by J.\,Shareshian
  and R.\,Woodroofe in \cite{SW}.  They proved that if $\TheGroup$ is a finite
  group, the coset poset $\CosPosetOf{\TheGroup}$ is not contractible.

  This motivated the second author of this paper to study the homotopy type of
  the coset poset of a finitely generated infinite group, focusing on the
  contractibility in her PhD thesis~\cite{Welsch}.  In her thesis, she was more
  interested in some special subests of the coset posets, than the coset poset
  itself, namely the finite index coset poset
  \(
    \CosPosetOf[\text{fi}]{\TheGroup} :=
    \CosPosetOf{\HerFamOf[\text{fi}]{\TheGroup}}
  \).
  She proved contractibility for many groups and provided nearly all
  examples of Section~\ref{one-point-two}.  Moreover the thesis showed that there
  are some groups with non-contractible finite index coset poset.  By her
  examples, she was led to conjecture Theorem~\ref{equiv}, which we are now able
  to prove using an extended version of her cone argument.

  Although Theorem~\ref{equiv} formally implies the result of
  Shareshian--Woodroofe that for finite $\TheGroup$ the coset nerve
  $\CosNerveOf{\TheGroup}$ is not contractible, our proof actually makes use of
  their result. In order to prove Theorem~\ref{normal-equiv}, we show the
  corresponding result for the family of proper normal subgroups, which might be
  of independent interest.
  \begin{NewTh}<statement>{Theorem~\ref{SW-analogon}}
    If $\TheGroup$ is finite, $\CosNerveOf[\text{nor}]{\TheGroup}$ is not
    contractible.
  \end{NewTh}

  \section{The family of finite index subgroups}
  In this section we characterize those groups $\TheGroup$, for which
  $\CosNerveOf[\text{fi}]{\TheGroup}$ is contractible. It will turn out that
  $\CosNerveOf[\text{fi}]{\TheGroup}$ is contractible if and only if the poset
  $\HerFamOf[\text{fi}]{\TheGroup}$ of proper finite index subgroups has
  infinitely many maximal elements. See Theorem~\ref{equiv} for the complete
  statement.

  \begin{observation}\label{cofinal}
    Let $\TheGrFam$ and $\AltGrFam$ be two \notion{co-final} families of
    subgroups in $\TheGroup$, i.e., for any subgroup $\TheSubgroup\in\TheGrFam$
    there is a subgroup $\AltSubgroup\in\AltGrFam$ containing $\TheSubgroup$, and
    vice versa.

    Then $\CosNerveOf{\TheGrFam}$ and $\CosNerveOf{\AltGrFam}$ are homotopy
    equivalent via the chain
    \[
      \CosNerveOf{\TheGrFam}
      \homotopic
      \CosSimplOf{\TheGrFam}
      =
      \CosSimplOf{\AltGrFam}
      \homotopic
      \CosNerveOf{\AltGrFam}
      \qed
    \]
  \end{observation}
  The family $\HerFamOf[\text{fi}][\text{max}]{\TheGroup}$ of maximal elements
  in $\HerFamOf[\text{fi}]{\TheGroup}$ is co-final in
  $\HerFamOf[\text{fi}]{\TheGroup}$, and for the corresponding coset nerve
  $\CosNerveOf[\text{max,fi}]{\TheGroup}$, we obtain the following
  homotopy equivalence:
  \begin{equation}\label{fin-hom}
    \CosNerveOf[\text{fi}]{\TheGroup}
    \homotopic
    \CosNerveOf[\text{max,fi}]{\TheGroup}
    .
  \end{equation}
  From this, we obtain one direction of our characterization.
  \begin{cor}\label{notcontr}
    Suppose that $\HerFamOf[\text{fi}]{\TheGroup}$ has only finitely many maximal
    finite index subgroups, i.e., the family
    \(
      \HerFamOf[\text{max,fi}]{\TheGroup}
      =
      \SetOf{
        \TheMaxSubgroup[1],\ldots,\TheMaxSubgroup[n]
      }
    \)
    is finite. Then, the Frattini subgroup
    \[
      \TheFratSubGr :=
      \TheMaxSubgroup[1] \intersection\cdots\intersection \TheMaxSubgroup[n]
    \]
    is a normal subgroup of finite index in $\TheGroup$ and the coset nerve
    $\CosNerveOf[\text{fi}]{\TheGroup}$ is homotopy equivalent to
    $\CosNerveOf{\TheGroup\rmod\TheFratSubGr}$.

    In particular, $\CosNerveOf[\text{max,fi}]{\TheGroup}$ is \emph{not} contractible.
  \end{cor}
  \begin{proof}
    The subgroups $\TheMaxSubgroup[i]$ are in $1$--$1$~correspondence to the
    maximal subgroups of the finite group $\TheGroup\rmod\TheFratSubGr$, and the
    projection $\TheGroup\rightarrow\TheGroup\rmod\TheFratSubGr$ induces an
    isomorphism of nerves
    \(
      \CosNerveOf[\text{max,fi}]{\TheGroup}
      \isomorphic
      \CosNerveOf[\text{max}]{\TheGroup\rmod\TheFratSubGr}
    \).
    The homotopy equivalence~\pref{fin-hom} implies
    \[
      \CosNerveOf[\text{fi}]{\TheGroup}
      \homotopic
      \CosNerveOf[\text{max,fi}]{\TheGroup}
      \isomorphic
      \CosNerveOf[\text{max}]{\TheGroup\rmod\TheFratSubGr}
      \homotopic
      \CosNerveOf{\TheGroup\rmod\TheFratSubGr}
    \]
    As $\TheGroup\rmod\TheFratSubGr$ is a nontrivial finite group, its coset
    nerve is not contractible as shown by Shareshian--Woodroofe~\cite{SW}.
  \end{proof}

  For the converse, we consider \notion{complementary subgroups}. We call a
  proper subgroup $\AltSubgroup\propersubgroup\TheGroup$ a \notion{complement} of
  the subgroup $\TheSubgroup\subgroup\TheGroup$ if
  $\TheGroup=\TheSubgroup\AltSubgroup$. Equivalently, one can say that
  $\AltSubgroup$ intersects each coset of $\TheSubgroup$.  Given a collection
  $\TheSubgroup[1],\ldots,\TheSubgroup[n]$, we say that $\AltSubgroup$ is a
  \notion{common complement} of the $\TheSubgroup[i]$ if it is a complement of
  their intersection
  \(
    \TheSubgroup[1]\intersect\cdots\intersect
    \TheSubgroup[n]
  \).
  Note that a complement $\AltSubgroup$ of a subgroup $\TheSubgroup$ cannot
  contain $\TheSubgroup$ because then
  $\TheGroup=\TheSubgroup\AltSubgroup=\AltSubgroup$ while we assume the complement
  $\AltSubgroup$ to be a proper subgroup of $\TheGroup$.

  We can now state the key argument of this paper.
  \begin{prop}[Cone construction]\label{cone}
    Let $\TheGrFam$ be a family of subgroups in $\TheGroup$ such that any
    finitely many subgroups $\TheSubgroup[1],\ldots,\TheSubgroup[n]\in\TheGrFam$ have
    a common complement $\AltSubgroup\in\TheGrFam$. Then the coset nerve
    $\CosNerveOf{\TheGrFam}$ is contractible.
  \end{prop}
  \begin{proof}
    Let $\TheSubComplex$ be a finite subcomplex of
    $\CosNerveOf{\TheGrFam}$. Thus, there are finitely many subgroups
    $\TheSubgroup[1],\ldots,\TheSubgroup[n]\in\TheGrFam$ such that
    \(
      \TheSubComplex
      \subcomplex
      \CosNerveOf{\SetOf{\TheSubgroup[1],\ldots,\TheSubgroup[n]}}
      \subcomplex
      \CosNerveOf{\TheGrFam}
    \).
    By hypothesis, there is a common complement $\AltSubgroup$ for the subgroups
    $\TheSubgroup[1],\ldots,\TheSubgroup[n]$. Then the coset $1\AltSubgroup$ is a
    vertex in $\CosNerveOf{\TheGrFam}$ whose star contains the whole subcomplex
    $\TheSubComplex$. Therefore $\TheSubComplex$ can be contracted within
    $\CosNerveOf{\TheGrFam}$.

    Since spheres are compact, any element of any homotopy group can be realized
    within a finite subcomplex of $\CosNerveOf{\TheGrFam}$. Thus, all homotopy
    groups of $\CosNerveOf{\TheGrFam}$ are trivial and $\CosNerveOf{\TheGrFam}$
    is contractible by Whitehead's Theorem.
  \end{proof}

  The main result of this section is an easy consequence.
  \begin{theorem}\label{equiv}
    For any group $\TheGroup$ the following are equivalent:
    \begin{enumerate}
      \item\label{A}
        Every proper subgroup of finite index in $\TheGroup$ has a complement
        that is also of finite index.
      \item\label{B}
        The nerve $\CosNerveOf[\text{fi}]{\TheGroup}$ is contractible.
      \item\label{C}
        The poset $\HerFamOf[\text{fi}]{\TheGroup}$ of proper finite index
        subgroups in $\TheGroup$ has infinitely many maximal elements.
    \end{enumerate}
  \end{theorem}
  \begin{proof}
    First, we show that \pref{A} implies \pref{B}. Let
    $\TheSubgroup[1],\ldots,\TheSubgroup[n]$ be proper finite index subgroups of
    $\TheGroup$. Their intersection
    \(
      \TheSubgroup[1]\intersect\cdots\intersect
      \TheSubgroup[n]
    \)
    is a proper finite index subgroup of $\TheGroup$, which by
    hypothesis~\pref{A} has a complement
    $\AltSubgroup\in\HerFamOf[\text{fi}]{\TheGroup}$. This is a common complement of
    the collection $\TheSubgroup[1],\ldots,\TheSubgroup[n]$. By the cone principle,
    we conclude that $\CosNerveOf[\text{fi}]{\TheGroup}$ is contractible.

    That \pref{B} implies \pref{C} has been proven in Corollary~\ref{notcontr}.

    For the remaining implication assume that
    $\HerFamOf[\text{fi}][\text{max}]{\TheGroup}$ is infinite. Let $\TheSubgroup$ be
    a proper subgroup of finite index in $\TheGroup$. We have to find a
    complementary subgroup $\AltSubgroup\in\HerFamOf[\text{fi}]{\TheGroup}$. Since a
    complement to a subgroup of $\TheSubgroup$ will also be a complement to
    $\TheSubgroup$, we may assume without loss of generality that $\TheSubgroup$ is
    normal in $\TheGroup$. Then, for any $\AltSubgroup\subgroup\TheGroup$, the
    product $\TheSubgroup\AltSubgroup$ is a subgroup of $\TheGroup$. Since only
    finitely many maximal subgroups contain the finite index subgroup $\TheSubgroup$
    (they correspond to the maximal subgroups of the finite group
    $\TheGroup\rmod\TheSubgroup$), we can choose a maximal subgroup $\AltSubgroup$
    not containing $\TheSubgroup$. Then $\TheGroup=\AltSubgroup\TheSubgroup$ since
    the maximal subgroup $\AltSubgroup$ is a proper subgroup of
    $\AltSubgroup\TheSubgroup$. Thus, we have found a complement for $\TheSubgroup$.
    Condition~\pref{A} follows.
  \end{proof}

  \subsection{Inheriting a contractible coset nerve}
  \begin{prop}\label{fin-index}
    Let $\TheGroup$ be a group with contractible coset nerve
    $\CosNerveOf[\text{fi}]{\TheGroup}$ and let $\TheSubgroup\subgroup\TheGroup$ be
    a subgroup of finite index. Then its coset nerve
    $\CosNerveOf[\text{fi}]{\TheSubgroup}$ is also contractible.
  \end{prop}
  \begin{proof}
    We use the characterization by the existence of complements.  Let
    $\TheSubSubgroup$ be a proper finite index subgroup of $\TheSubgroup$. Then it
    is a proper finite index subgroup of $\TheGroup$. Since
    $\CosNerveOf[\text{fi}]{\TheGroup}$ is contractible, there is a finite
    index complement $\AltSubgroup$ for $\TheSubSubgroup$ in $\TheGroup$. Then
    $\AltSubgroup\intersect\TheSubgroup$ is a complement for $\TheSubSubgroup$ in
    $\TheSubgroup$.
  \end{proof}
  \begin{observation}\label{projektion}
    Let $\TheProj\mapcolon\AltGroup\rightarrow\TheGroup$ be an epimorphism of groups
    and let $\TheSubgroup\subgroup\TheGroup$ be a maximal subgroup. Then, the preimage
    $\TheProjOf[][-1]{\TheSubgroup}$ is a maximal subgroup of $\AltGroup$.

    Hence, $\CosNerveOf[\text{fi}]{\AltGroup}$ is contractible provided that
    $\CosNerveOf[\text{fi}]{\TheGroup}$ is contractible.\qed
  \end{observation}

  \subsection{Examples: groups with contractible coset nerve}\label{one-point-two}
  A group $\TheGroup$ is called \notion{indicable} if it admits an epimorphism
  onto the infinite cyclic group $\IntNumbers$. Since $\IntNumbers$ has
  infinitely many maximal subgroups, its coset nerve is contractible; and $\TheGroup$
  inherits this property by Observation~\ref{projektion}.

  Hence the following groups all have contractible coset nerves:

  \begin{enumerate}
    \item
      Free groups.
    \item
      $\TheGroup\freeprod\IntNumbers$,
      $\TheGroup\crossprod\IntNumbers$,
      $\TheGroup\semidirectprod\IntNumbers$ for any group $\TheGroup$.
    \item
      Free abelian groups.
    \item
      HNN extensions.
    \item
      Orientation preserving Fuchsian groups of genus at least~1. This can be seen
      by abelianizing the presentation.
    \item
      Other Fuchsian groups, provided the genus is at least $2$.
    \item
      Baumslag--Solitar groups $\GroupPresented[m,n]{ba^m=a^nb}$.
    \item
      Thompson's group $F$ has $\IntNumbers\crossprod\IntNumbers$ as its maximal
      abelian quotient.
    \item
      Artin groups.
    \item
      Pure braid groups (finite index subgroups of indicable groups are indicable).
  \end{enumerate}

  As a variation of this trick, we can use the infinite dihedral group $\DiederGr[\infty]$
  instead of $\IntNumbers$. Thus, all of the following groups also have a contractible
  coset nerve:
  \begin{enumerate}[resume]
    \item
      Semi-direct products $\IntNumbers\semidirectprod\TheGroup$ map
      to $\DiederGr[\infty]$ with full image or image $\IntNumbers$.
    \item
      Infinite Coxeter groups
      \[
      \GroupPresented[ s_1,\ldots,s_n ]{
      s_i^2,\,
      (s_i s_j)^{m_{ij}}
      },
      \]
      where all $m_{ij}$ are even and $m_{1,2}=\infty$. Killing all $s_i$ for
      $i \geq 3$, we recognize $\DiederGr[\infty]$ as a quotient.
  \end{enumerate}

  We can generalize from the infinite dihedral group to other euclidean
  crystallographic groups.
  \begin{observation}
    Let $\TheGroup$ be a crystallographic group in $\RealNumbers[][\TheDim]$.
    Then $\TheGroup$ has finite index in a semi-direct product
    \(
      \AltGroup=\IntNumbers[][\TheDim]
      \semidirectprod
      \FinGroup
    \)
    where the finite quotient $\FinGroup$ acts on the lattice $\IntNumbers[][\TheDim]$
    via linear automorphisms. Hence the rescaled lattices
    $\TheSubgroup[\TheNumber]:=\TheNumber\IntNumbers[][\TheDim]$
    are invariant sublattices. Thus, the group
    \(
      \TheSubgroup[\TheNumber]
      \semidirectprod
      \FinGroup
    \)
    is a subgroup of index $\TheNumber[][\TheDim]$ in $\AltGroup$. In particular,
    for each prime number $p$, the group $\AltGroup$ has a subgroup of index
    $p^{\TheDim}$. Thus, the index of a maximal subgroup containing
    \(
      \TheSubgroup[\TheNumber]
      \semidirectprod
      \FinGroup
    \)
    is a $p$-power. It follows that $\AltGroup$ has infinitely many maximal
    subgroups.

    By Theorem~\ref{equiv}, the coset nerve of $\AltGroup$ is contractible.
    By Proposition~\ref{fin-index}, this carries over to the coset nerve of
    the finite index subgroup $\TheGroup\subgroup\AltGroup$.
  \end{observation}

  Dealing with $\SlOf[n]{\IntNumbers}$ requires a new idea.
  \begin{observation}\label{einfach}
    Let $\TheGroup$ be a finitely generated group that has infinitely many
    finite index normal subgroups with simple quotient. Then any proper
    finite index subgroup $\TheSubgroup\subgroup\TheGroup$ has a finite index
    complement. In particular, the coset nerve $\CosNerveOf[\text{fi}]{\TheGroup}$
    is contractible.
  \end{observation}
  \begin{proof}
    $\TheGroup$ acts by left multiplication on the finite quotient
    $\TheGroup\rmod\TheSubgroup$. Let $\TheKernel$ be the kernel of this
    action. Thus, $\TheKernel$ is a normal subgroup of finite index
    in $\TheGroup$. For any normal subgroup $\TheNormalSubGr$ of $\TheGroup$,
    the product $\TheKernel\TheNormalSubGr$ is normal in $\TheGroup$.
    If the quotient $\TheGroup\rmod\TheNormalSubGr$ is simple and
    if $\TheNormalSubGr$ does not contain $\TheKernel$, then
    $\TheGroup=\TheKernel\TheNormalSubGr=\TheSubgroup\TheNormalSubGr$.

    As $\TheKernel$ has finite index, it is contained only within finitely many
    normal subgroups $\TheNormalSubGr$ of $\TheGroup$. Thus, there are (infinitely
    many) normal subgroups $\TheNormalSubGr$ that are complementary to $\TheKernel$
    and hence to $\TheSubgroup$.
  \end{proof}
  Examples of groups whose coset nerve is contractible by this argument include:
  \begin{enumerate}[resume]
    \item
      $\SlOf[n]{\IntNumbers}$,
      $\AutOf{\FreeGroup[n]}$,
      and
      $\OutOf{\FreeGroup[n]}$.
    \item
      Symplectic groups and mapping class groups.
  \end{enumerate}

  \subsection{Examples: groups with non-contractible coset nerves}
  As infinite simple groups do not have any proper subgroups of finite index,
  they can be used to construct some easy examples of groups with non-contractible
  coset nerve.
  \begin{observation}
    Suppose the infinite group $\TheGroup$ has a simple subgroup
    $\TheSimpleGr$. Let $\TheNormalSubGr$ be the normal subgroup generated by
    $\TheSimpleGr$.  Then $\TheGroup$ and $\TheGroup\rmod\TheNormalSubGr$ have the
    same finite quotients since a projection homomorphism
    $\TheProj\mapcolon\TheGroup\rightarrow\FinGroup$ is trivial on $\TheSimpleGr$
    (by simplicity) and hence is trivial on $\TheNormalSubGr$.

    Hence $\TheGroup$ has only finitely many normal subgroups of finite index
    provided that $\TheGroup\rmod\TheNormalSubGr$ has only finitely many normal
    subgroups of finite index. As each such normal subgroup is contained in only
    finitely many subgroups of $\TheGroup$, the whole family
    $\HerFamOf[\text{fi}]{\TheGroup}$ is finite in this case; and the coset nerve
    $\CosNerveOf[\text{fi}]{\TheGroup}$ is not contractible.
  \end{observation}
  This applies for instance to any group of the form
  \[
    \text{infinite simple}
    \freeprod
    ( \text{finite} \crossprod \text{infinite simple} )
  \]
  Of course, this theme allows for many variations.
  More interesting are residually finite examples. The examples we know are
  $p$-groups.
  \begin{observation}
    Let $p$ be a prime. Any finitely generated $p$-group $\TheGroup$ has only
    finitely many maximal subgroups. In particular, $\CosNerveOf[\text{fi}]{\TheGroup}$
    is not contractible.
  \end{observation}
  \begin{proof}
    We claim that every maximal finite index subgroup in $\TheGroup$ has index $p$.
    Since $\TheGroup$ is finitely generated, there are only finitely many subgroups
    of index $p$.

    A maximal finite index subgroup $\TheSubgroup\subgroup\TheGroup$ contains
    a finite index subgroup $\TheNormalSubGr$ which is normal in $\TheGroup$.
    Hence the maximal subgroup $\TheSubgroup$ corresponds to a maximal subgroup
    of the finite $p$-group $\TheGroup\rmod\TheNormalSubGr$. Therefore, $\TheSubgroup$
    has index $p$ in $\TheGroup$.
  \end{proof}
  Thus, the first Grigorchuk group or the Gupta--Sidki groups have non-contractible
  coset nerves.

  \section{The family of finite index normal subgroups}
  In this section, we consider the family
  \[
    \HerFamOf[\text{nor,fi}]{\TheGroup}
    :=
    \SetOf[{
      \TheNormalSubGr\propernormalsubgr\TheGroup
    }]{
      \TheNormalSubGr\text{\ has finite index in\ }\TheGroup
    }
  \]
  of proper normal subgroups that have finite index in $\TheGroup$. Let
  $\CosNerveOf[\text{nor,fi}]{\TheGroup}$ be the coset nerve for this family.  We
  shall prove the exact analogue of Theorem~\ref{equiv}. For this we need the
  analogue of the result of Shareshian--Woodroofe~\cite{SW}, which we shall prove
  as Theorem~\ref{SW-analogon} in the appendix.

  \begin{theorem}\label{normal-equiv}
    For any group $\TheGroup$, the following are equivalent:
    \begin{enumerate}
      \item\label{a}
        $\CosNerveOf[\text{nor,fi}]{\TheGroup}$ is contractible.
      \item\label{b}
        $\HerFamOf[\text{nor,fi}]{\TheGroup}$ has infinitely many maximal elements.
      \item\label{c}
        Any proper finite index subgroup $\TheSubgroup\subgroup\TheGroup$ has
        infinitely many complementary normal subgroups of finite index
        in $\TheGroup$.
      \item\label{d}
        Each proper normal subgroup of finite index in $\TheGroup$ has a complementary
        normal subgroup of finite index in $\TheGroup$.
    \end{enumerate}
  \end{theorem}
  \begin{proof}
    To see that \pref{a} implies \pref{b}, we argue the contrapositive. So
    assume that $\TheGroup$ has only finitely many maximal proper normal subgroups
    $\MaxNormalSubGr[1],\ldots,\MaxNormalSubGr[n]$ of finite index. Their
    intersection
    \[
      \NorFratSubGr :=
      \MaxNormalSubGr[1]\intersect\cdots\intersect\MaxNormalSubGr[n]
    \]
    is normal and of finite index. By correspondence, the coset nerves
    $\CosNerveOf[\text{nor,fi}]{\TheGroup}$ and
    $\CosNerveOf[\text{nor,fi}]{\TheGroup\rmod\NorFratSubGr}$ are isomorphic.
    By Theorem~\ref{SW-analogon} from the appendix, the latter is not contractible.

    Now we assume \pref{b} and show \pref{c}. So let
    $\TheSubgroup\subgroup\TheGroup$ be a proper subgroup of finite index. As
    before, we consider the action of $\TheGroup$ on the finite set
    $\TheGroup\rmod\TheSubgroup$. The kernel $\TheKernel$ of the action is a proper
    normal subgroup of finite index in $\TheGroup$. Any maximal proper normal
    subgroup $\TheNormalSubGr\propernormalsubgr\TheGroup$ not containing
    $\TheKernel$ is a complement to $\TheKernel$ and thus to $\TheSubgroup$.  Since
    $\TheKernel$ is only contained in finitely many such $\TheNormalSubGr$, there
    are infinitely many complementary subgroups in
    $\HerFamOf[\text{nor,fi}]{\TheGroup}$.

    Clearly \pref{c} implies \pref{d}.

    The remaining implication from \pref{d} to \pref{a} follows immediately from
    the cone construction~\ref{cone} applied to the family
    $\HerFamOf[\text{nor,fi}]{\TheGroup}$.
  \end{proof}
  Under the equivalent conditions of Theorem~\ref{normal-equiv} each
  proper finite index subgroup has a complement.  Thus, we obtain the following:
  \begin{cor}
    If $\CosNerveOf[\text{nor,fi}]{\TheGroup}$ is contractible, then so is
    $\CosNerveOf[\text{fi}]{\TheGroup}$.
  \end{cor}

  For finitely generated $\TheGroup$, we also have a characterization via finite
  simple quotients.
  \begin{observation}
    If $\TheGroup$ is finitely generated, the following are equivalent:
    \begin{enumerate}
      \item\label{e}
        $\HerFamOf[\text{nor,fi}]{\TheGroup}$ has infinitely many maximal
        elements.
      \item\label{f}
        $\TheGroup$ has infinitely many pairwise non-isomorphic finite simple
        quotients.
    \end{enumerate}
  \end{observation}
  \begin{proof}
    The direction from \pref{f} to \pref{e} is obvious: for each finite simple
    quotient, the kernel of the projection is a maximal finite index normal
    subgroup.

    For the converse, we just observe that for a finite (simple) group
    $\FinGroup$, there are only finitely many homomorphisms
    \(
      \TheGroup\rightarrow\FinGroup
    \)
    as $\TheGroup$ is finitely generated. Hence, each such possible quotient
    arises from at most finitely many normal subgroups.
  \end{proof}

  \section{The family of all proper subgroups}\label{all-subgroups}
  In the poset of proper finite index subgroups, the maximal elements form
  a co-final family. Our analysis in that case crucially relied on that
  feature. For finitely generated groups $\TheGroup$, the poset of all
  proper subgroups behaves in the same way.
  \begin{lemma}\label{higher-finite-index}
    Let $\TheGroup$ be a finitely generated group. Then any proper subgroup
    of $\TheGroup$ is contained in a maximal proper subgroup. In particular,
    the families $\HerFamOf{\TheGroup}$ and $\HerFamOf[][\text{max}]{\TheGroup}$
    are co-final, and we have the corresponding homotopy equivalence:
    \[
      \CosNerveOf{\TheGroup} \homotopic  \CosNerveOf[\text{max}]{\TheGroup}
    \]
  \end{lemma}
  \begin{proof}
    By the Lemma of Zorn, it suffices to see that
    any ascending union of proper subgroups in $\TheGroup$ is a proper subgroup.
    However, this follows from the existence of a finite generating set for
    $\TheGroup$. If an ascending union contained all the generators, this would
    happen already at some stage along the underlying ascending chain of subgroups.
  \end{proof}
  \begin{observation}
    If a maximal subgroup
    $\TheMaxSubgroup\subgroup\TheGroup$ has only finitely many conjugate subgroups in
    $\TheGroup$, it is of finite index in $\TheGroup$.
  \end{observation}
  \begin{proof}
    The normalizer of $\TheNormalizerOf{\TheMaxSubgroup}$ of $\TheMaxSubgroup$
    has finite index in $\TheGroup$ as it is a stabilizer of a $\TheGroup$-action on
    a finite set. Moreover, we have the inclusions:
    \[
      \TheMaxSubgroup \subgroup
      \TheNormalizerOf{\TheMaxSubgroup} \subgroup \TheGroup
    \]
    As $\TheMaxSubgroup$ is maximal, it follows that
    $\TheMaxSubgroup=\TheNormalizerOf{\TheMaxSubgroup}$ (and therefore has finite index in
    $\TheGroup$) or that $\TheMaxSubgroup$ is normal in $\TheGroup$.

    If $\TheMaxSubgroup$ is normal, it follows that the trivial group is a maximal
    subgroup in the quotient $\TheGroup\rmod\TheMaxSubgroup$. Hence the quotient is
    cyclic of prime order. Again $\TheMaxSubgroup$ has finite index in $\TheGroup$.
  \end{proof}
  Since the collection of maximal subgroups is closed with respect to
  taking conjugates, we infer the following:
  \begin{cor}\label{endlich-viele-dann-alle-koendlich}
    If a finitely generated group $\TheGroup$ has only finitely many maximal
    subgroups, they all are of finite index in $\TheGroup$ and the coset nerve
    $\CosNerveOf{\TheGroup}$ is homotopy equivalent to $\CosNerveOf[\text{fi}]{\TheGroup}$.
    Both are non-contractible.\qed
  \end{cor}

  We conclude this section with two results towards contractibility.
  \begin{observation}
    If every proper subgroup of $\TheGroup$ is contained in a proper finite
    index subgroup of $\TheGroup$, then $\HerFamOf{\TheGroup}$ and
    $\HerFamOf[\text{fi}]{\TheGroup}$ are co-final. In this case, we have the
    homotopy equivalence
    \[
      \CosNerveOf{\TheGroup}
      \homotopic
      \CosNerveOf[\text{fi}]{\TheGroup}
    \]
    and contractibility of $\CosNerveOf[\text{fi}]{\TheGroup}$ implies
    contractibility of $\CosNerveOf{\TheGroup}$.\qed
  \end{observation}
  A slightly more involved argument is required in the following situation.
  We call a subgroup $\TheSubgroup$ of $\TheGroup$ \notion{extreme}
  if it is of finite index or not even contained in a subgroup of finite
  index. So, extreme subgroups are either very large or very small.
  \begin{prop}
    Let $\TheGroup$ be a group and assume that every finite collection
    $\TheSubgroup[1],\ldots,\TheSubgroup[n]\subgroup\TheGroup$ of proper
    extreme subgroups has a common complement.
    Then the coset nerve $\CosNerveOf{\TheGroup}$ is contractible.
  \end{prop}
  \begin{proof}
    We consider the family $\TheGrFam$ of proper extreme subgroups of
    $\TheGroup$. This family is obtained from $\HerFamOf{\TheGroup}$ by removing
    all subgroups of infinite index that are contained in a proper finite index
    subgroup of $\TheGroup$. Clearly, $\TheGrFam$ and $\HerFamOf{\TheGroup}$ are
    co-final, wherefore
    \(
      \CosNerveOf{\TheGrFam}\homotopic\CosNerveOf{\TheGroup}
    \).

    By the cone construction~\ref{cone}, it suffices to find for any finite
    collection $\TheSubgroup[1],\ldots,\TheSubgroup[n]\in\TheGrFam$ a common
    complement in $\TheGrFam$. By hypothesis, the $\TheSubgroup[i]$ have
    a common complementary subgroup $\AltSubgroup\subgroup\TheGroup$.
    If $\AltSubgroup\in\TheGrFam$, there is nothing to be argued.
    Otherwise, $\AltSubgroup$ has infinite index in $\TheGroup$ and is
    contained in a proper finite index subgroup, which then can be used
    as a common complement for the $\TheSubgroup[i]$.
  \end{proof}

  \section{Appendix: the coset nerve for proper normal subgroups in a finite group}
  In this appendix, we consider the family
  \[
    \HerFamOf[\text{nor}]{\TheGroup}
    :=
    \SetOf[{
      \TheNormalSubGr\normalsubgr\TheGroup
    }]{
      \TheNormalSubGr\neq\TheGroup
    }
  \]
  of proper normal subgroups of $\TheGroup$. Let $\CosNerveOf[\text{nor}]{\TheGroup}$
  denote the coset nerve associate to $\HerFamOf[\text{nor}]{\TheGroup}$. This section
  is devoted to a proof of the
  following analogue to the result of Shareshian--Woodroofe~\cite{SW}:
  \begin{theorem}\label{SW-analogon}
    If $\TheGroup$ is finite, $\CosNerveOf[\text{nor}]{\TheGroup}$ is not
    contractible.
  \end{theorem}
  \begin{proof}
    If the group $\TheGroup$ is trivial, the family of proper normal subgroups
    is empty. Hence, we assume that $\TheGroup$ is nontrivial.

    The family $\HerFamOf[\text{nor}]{\TheGroup}$ is closed with respect to
    intersections. Hence, we shall consider the associated order complex
    $\CosPosetOf[\text{nor}]{\TheGroup}$ instead of the homotopy equivalent
    coset nerve. Our argument uses many ideas of Brown~\cite[Section~8]{Brown}.

    Inducting on the complexity of $\TheGroup$, we shall show that
    $\CosPosetOf[\text{nor}]{\TheGroup}$ is a nontrivial wedge of spheres.  Simple
    groups make up the base of the induction. Fortunately, simple groups are easily
    understood since the only proper normal subgroup is the trivial subgroup.
    \begin{claim}\label{claim-one}
      If $\TheGroup$ is simple, $\CosPosetOf[\text{nor}]{\TheGroup}$ is a finite
      discrete set, which we shall identify with $\TheGroup$ itself.
      It is $0$-spherical and not contractible (since $\TheGroup$ is
      nontrivial).
    \end{claim}

    Now assume that $\TheGroup$ is not simple. Let
    $\MinNormalSubGr\propernormalsubgr\TheGroup$ be a minimal normal subgroup.
    Let
    \[
      \TheProjection \mapcolon \TheGroup \longrightarrow
      \TheGroup\rmod \MinNormalSubGr
    \]
    denote the canonical projection.
    \begin{claim}\label{claim-two}
      Assume that $\TheGroup$ contains no proper normal subgroup that surjects
      onto the quotient $\TheGroup\rmod\MinNormalSubGr$. Then the map
      \begin{align*}
        \TheNerveMap \mapcolon \CosPosetOf[\text{nor}]{\TheGroup}
        & \longrightarrow \CosPosetOf[\text{nor}]{\TheGroup\rmod\MinNormalSubGr} \\
        \TheGroupElement \TheNormalSubGr & \mapsto
                                           \TheProjectionOf{ \TheGroupElement \TheNormalSubGr }
      \end{align*}
      is a homotopy equivalence.

      Thus sphericity and non-contractibility of $\CosPosetOf[\text{nor}]{\TheGroup}$
      is inherited from $\CosPosetOf[\text{nor}]{\TheGroup\rmod\MinNormalSubGr}$.
    \end{claim}
    To see this, consider the intersection-closed family
    \(
      \TheGrFam := \SetOf[{
        \TheNormalSubGr \propernormalsubgr \TheGroup
      }]{
        \MinNormalSubGr \subseteq \TheNormalSubGr
      }
    \).
    The elements of $\TheGrFam$ are in $1$--$1$ correspondence to the proper normal
    subgroups of $\TheGroup\rmod\MinNormalSubGr$. Therefore,
    \(
      \CosPosetOf[\TheGrFam]{\TheGroup}
    \)
    and
    \(
      \CosPosetOf[\text{nor}]{\TheGroup\rmod\MinNormalSubGr}
    \)
    are isomorphic.  On the other hand, $\TheGrFam$ is co-final in
    $\HerFamOf[\text{nor}]{\TheGroup}$ since by hypothesis
    $\TheNormalSubGr\MinNormalSubGr$ is a proper normal subgroup of $\TheGroup$ for
    any proper normal $\TheNormalSubGr\propernormalsubgr\TheGroup$.
    This proves Claim~\ref{claim-two}.

    It remains to deal with the case that there is a proper normal subgroup
    $\FixNormalSubGr\propernormalsubgr\TheGroup$ that surjects onto
    $\TheGroup\rmod\MinNormalSubGr$. Note that the intersection
    $\FixNormalSubGr\intersect\MinNormalSubGr$ is normal in $\TheGroup$. As
    $\FixNormalSubGr$ is proper, but the product $\FixNormalSubGr\MinNormalSubGr$ is
    all of $\TheGroup$, we conclude that
    $\FixNormalSubGr\intersect\MinNormalSubGr\neq \MinNormalSubGr$. As
    $\MinNormalSubGr$ is a minimal normal subgroup in $\TheGroup$, we find:
    \[
      \FixNormalSubGr\intersect\MinNormalSubGr = \SetOf{1}
      \qquad\text{and}\qquad
      \TheGroup=\MinNormalSubGr\crossprod\FixNormalSubGr
    \]
    In particular, we can identify
    $\FixNormalSubGr$ with $\TheGroup\rmod\MinNormalSubGr$. We call a coset
    $\TheGroupElement\TheSubgroup$ in $\TheGroup$ \notion{large} if it surjectsq
    onto $\FixNormalSubGr$. We call it \notion{small} otherwise.
    Note that the coset $\TheGroupElement\TheSubgroup$ is small if and only if
    $1\TheSubgroup$ is small.

    We consider the following families:
    \begin{align*}
      \HerFamOf[\text{large}]{\TheGroup}
      & := \SetOf[ \TheNormalSubGr \propernormalsubgr \TheGroup]{
        1\TheNormalSubGr \text{\ is large\,}
        }
      \\
      \HerFamOf[\text{small}]{\TheGroup}
      & := \SetOf[ \TheNormalSubGr \propernormalsubgr \TheGroup]{
        1\TheNormalSubGr \text{\ is small\,}
        }
      \\
      \PairFam
      & := \SetOf[{
        \MinNormalSubGr \crossprod L
        }]{
        L \propernormalsubgr \FixNormalSubGr
        }
    \end{align*}
    Note that $\HerFamOf[\text{small}]{\TheGroup}$ and $\PairFam$ are
    closed with respect to intersections. Moreover
    $\PairFam$ is co-final in $\HerFamOf[\text{small}]{\TheGroup}$.
    More precisely for $\TheNormalSubGr\in\HerFamOf[\text{small}]{\TheGroup}$,
    we have
    \(
      \TheNormalSubGr
      \subseteq
      \MinNormalSubGr
      \crossprod
      \TheProjectionOf[\FixNormalSubGr]{\TheNormalSubGr}
      \in\PairFam
    \).
    Hence, we have the following homotopy equivalence:
    \begin{equation}\label{xxx}
      \CosPosetOf{
        \HerFamOf[\text{small}]{\TheGroup}
      }
      \homotopic
      \CosPosetOf{\PairFam}
      \isom
      \CosPosetOf[\text{nor}]{\FixNormalSubGr}
    \end{equation}

    Now consider a large normal subgroup
    $\LargeNormalSubGr\propernormalsubgr\TheGroup$.  We claim that the projection
    onto $\FixNormalSubGr$ restricts to an isomorphism on $\LargeNormalSubGr$. The
    reason is that $\MinNormalSubGr$ is minimal and because of this
    $\LargeNormalSubGr\intersect\MinNormalSubGr$ is trivial or all of
    $\MinNormalSubGr$. However, the latter possibility is excluded since
    $\LargeNormalSubGr$ is a proper subgroup of $\TheGroup$. This argument shows:
    \begin{claim}
      If $\LargeNormalSubGr$ is a large proper normal subgroup of $\TheGroup$, then
      $\TheGroup=\MinNormalSubGr\crossprod\LargeNormalSubGr$.
    \end{claim}
    Consequently, all large proper normal subgroups have the same cardinality.
    Thus they are mutually incomparable with respect to inclusion.

    We now see the structure of the coset poset $\CosPosetOf[\text{nor}]{\TheGroup}$.
    The vertices from $\HerFamOf[\text{small}]{\TheGroup}$ form the \notion{bottom part}
    whose
    order complex is homotopy equivalent to
    \(
      \CosPosetOf[\text{nor}]{\FixNormalSubGr}
    \).
    Even more is true. The coset $1\FixNormalSubGr$ is large. Its link is spanned
    by small cosets contained in $1\FixNormalSubGr$. Therefore, this link is
    isomorphic to
    \(
      \CosPosetOf[\text{nor}]{\FixNormalSubGr}
    \)
    and the above homotopy equivalence~\ref{xxx} is a deformation retraction
    of the bottom part
    \(
      \CosPosetOf{
        \HerFamOf[\text{small}]{\TheGroup}
      }
    \)
    onto this link.

    The vertices from $\HerFamOf[\text{large}]{\TheGroup}$ lie above the bottom
    part. The order complex $\CosPosetOf[\text{nor}]{\TheGroup}$ is the union of the
    bottom part and the stars of the large vertices. Each such star is just the cone
    over the link of the vertex; and we may consider these stars independently
    because two large vertices are never joined by an edge.  Moreover, the link of a
    large coset $\TheGroupElement\LargeNormalSubGr$ is its \notion{descending link},
    i.e., all vertices in the link are cosets that are contained in
    $\TheGroupElement\LargeNormalSubGr$.
    \begin{claim}
      The link of a large vertex $\TheGroupElement\LargeNormalSubGr$ is
      isomorphic to $\CosPosetOf[\text{nor}]{\LargeNormalSubGr}$.
    \end{claim}
    This just follows from the fact that a coset
    $\TheGroupElement\TheNormalSubGr$ is contained in $\TheGroupElement\LargeNormalSubGr$
    if and only if $\TheNormalSubGr\subgroup\LargeNormalSubGr$.

    Now, we build $\CosPosetOf[\text{nor}]{\TheGroup}$ by adding the stars of
    large vertices, one by one, to the bottom part
    \(
      \CosPosetOf{
        \HerFamOf[\text{small}]{\TheGroup}
      }
    \).
    Adding $1\FixNormalSubGr$ as the first large vertex, we cone of its
    descending link. At this point, we obtain a contractible space.  Adding each of
    the remaining large vertices (and there are at least the other cosets of
    $\FixNormalSubGr$) amounts to wedging on the suspension of its link, i.e.,
    wedging on a copy of
    \(
      \SuspensionOf{
        \CosPosetOf[\text{nor}]{\FixNormalSubGr}
      }
    \).
    Thus, we have argued the following:
    \begin{claim}\label{claim-three}
      Assume that there is a proper normal subgroup
      $\FixNormalSubGr\propernormalsubgr\TheGroup$ that surjects onto
      $\TheGroup\rmod\MinNormalSubGr$. Then $\CosPosetOf[\text{nor}]{\TheGroup}$ is
      homotopy equivalent to a nontrivial wedge of copies of the suspension
      \(
        \SuspensionOf{
          \CosPosetOf[\text{nor}]{\FixNormalSubGr}
        }
      \).
    \end{claim}

    Note that the theorem follows by induction from the Claim~\ref{claim-one},
    Claim~\ref{claim-two}, and Claim~\ref{claim-three}.
  \end{proof}

  \medskip

  \paragraph{Acknowledgments.} We would like to thank Russ~Woodroofe for
  suggesting that a proof of Theorem~\ref{SW-analogon} could be carried out
  using ideas from Brown~\cite{Brown} rather than following the much steeper path
  of Shareshian--Woodroofe~\cite{SW}. We also thank Benjamin~Brück, Linus~Kramer,
  and Russ~Woodroofe for helpful comments on a preliminary version of this paper.

\end{document}